\documentclass{article}
\usepackage{amsmath}
\usepackage{hyperref}
\usepackage {amsthm}
\usepackage{amssymb}
\usepackage{amsfonts}
\usepackage{mathrsfs}
\usepackage {graphicx,color}
\theoremstyle{definition}
\newtheorem{lemma}{Lemma}[section]
\newtheorem{definition}[lemma]{Definition}
\newtheorem{proposition}[lemma]{Proposition}
\newtheorem{theorem}[lemma]{Theorem}

\newtheorem{example}{Example}
\newtheorem{remark}{Remark}
\newtheorem{question}{Question}
\newtheorem{claim}{Claim}

\numberwithin{equation}{section}
\renewcommand{\proof}{Proof. }
\renewcommand{\qed}{\hfill\small{$\square$}\normalsize}
\title{Convergence of metrics under self - dual Weyl tensor and Scalar curvature bounds}
\author{Yiyan Xu}
\begin{document}
\maketitle
\begin{abstract}
We establish a $C^{1,\alpha}$ compactness theorem for the metrics
with bounded self - dual Weyl tensor and Scalar curvature. The key
step is to estimate the $C^{1,\alpha}$ harmonic radius, where we use
the blow up analysis as in \cite{Anderson90}. The result is
motivated by, and may be applied to the Calabi flow on complex
surfaces \cite{Calabi82}.
\end{abstract}
\section{Introduction}
In Riemannian Geometry, it is interesting to consider the
convergence and finiteness in the presence of basic geometric data.
With the fundamental Cheeger - Gromov convergence theory, this
problem is reduced to the harmonic radius estimate in a standard
way.  In a remarkable piece of work, Jost and Karcher obtained a
explicit estimate on the harmonic radius which depends only on lower
volume, upper diameter and sectional curvature bound \cite{Jost84}
\cite{GW88}. The most important feature about harmonic coordinates
is that the metric is apparently controlled by the Ricci curvature,
roughly speaking, to a determined system of P.D.E., whereas assuming
bounds on $Rm$ corresponds to a overdetermined system information.
This was exploited by \cite{Anderson90}, i.e. the bounds $|Ric|\leq
\lambda$ and $\hbox{inj}\geq i_0$ alone imply a lower bound on the
harmonic radius without any assumptions on the curvature tensor.

Without the lower bound on the injectivity radius but relaxes to a
lower bound on volume, there are examples which exhibit that the
degeneration may occur. Pursuing this phenomena of thought,
geometers had made a well known conjecture: the convergence is
$C^{1,\alpha}$ off a singular set of Hausdorff codimension at least
four \cite{CCT02}.

Key steps in Andersen's proof are to establish the continuous
properties of the harmonic radius, regard as a function, under the
$C^{k,\alpha}(k\geq1)$ or $W^{k,p}(k\geq 1)$ convergence, which in
turn ultimately depends on some type of partial differential
equations estimate. Besides the elliptic regularity theory, the main
ingredient in the Andersen's results is a blow up argument that
relates the validity of convergence/rigidity theorems to the global
behavior of complete, non - compact manifolds with critical metric.
The hypothesis on the manifolds, such as curvature bound and lower
bound on injectivity radius, was basically only used in the
characterization of the Euclidean space, i.e. rigidity theorem. It
is likely that some weakening of the hypothesis, but still valid in
the rigidity theorem, suffices to give a lower bound on the harmonic
radius.

We proceed with another generalization of the metrics with bounded
Ricci curvature to scalar curvature bound, but at the expense of
some hypotheses on the structure of metric, for example,  Anti -
self - dual or  K\"{a}hler metric.  The motivation here also comes
from the Calabi flow. In estimating whether the solution of Calabi
flow exists and convergence, one crucial estimate which still
lacking is one showing that the space of all k\"{a}hler metrics
$\omega$ in a fixed cohomoplogy class $\Omega$, and for which the
scalar curvature $|S(\omega)|_{\omega}$ is bounded, is compact in
the $C^{1,\alpha}$ Cheeger - Gromov topology, see \cite{Calabi82}
and \cite{ChHe08}.

\theorem \label{MTCC} Let $M^4$ be a  compact oriented four
dimensional manifold, let $\{g_i\}$ be a sequence metrics on $M^4$
with
\begin{enumerate}
  \item either bounded self - dual Weyl tensor, i.e., $|W^{+}(g_i)|\leq
  \Lambda$ or $(M, g_i)$ is K\"{a}hler metric;
  \item bounded scalar curvature, i.e., $|S(g_i)|\leq \Lambda;$
  \item unit volume, i.e., $|\textrm{Vol}(M_i,g_i)|\equiv 1$;
  \item bounded Sobolev constant  $C_S(g_i)\leq C_S$, i.e.,
\begin{equation}\label{USobIM}
\{\int_M|v|^4d\nu_{g_i}\}^{\frac{1}{2}}\leq C_S\int_M|d
v|_{g_i}^2d\nu_{g_i}, \forall v\in C_0^{0,1}(M).
\end{equation}
\end{enumerate}
Then there exist a subsequence $\{j\}\subset\{i\}$  such that
$(M,g_j)$ converge to  a compact multi - fold $(M_\infty,g_\infty)$
in the Gromov - Hausdorff topology. Moreover, on the regular set
$M_\infty\setminus\{x_1,\cdots,x_m\}$, the metric $g_\infty$ is
$C^{1,\alpha}$ and the convergence is in the $C^{1,\alpha}$ Cheeger
- Gromov topology; while each of the singular point $x_i$ have a
neighborhood homeomorphic to finite disjoint union of cones
$C(S^3/\Gamma)$ with identifications of vertex, where $\Gamma$ is
finite subgroup of $O(4)$.

\remark By the Gauss - Bonnet formula and the Hirzebruch signature
formula \cite{Besse87}, then we have
\begin{equation}\label{RmEnTop}
\int_{M^4}|Rm|^2=-24\pi^2\tau(M)-8\pi^2\chi(M)+\int_M\frac{1}{12}S^2+4|W^+|^2
\end{equation}
If the scalar curvature $S$ and self - dual Weyl tensor has
uniformly bounded $L^2-$ norm, then we have a prior $L^2-$ bound on
curvature tensor $Rm$. The above theorem holds, with out fixing
topology, under the additional assumption that
\[\int_{M_i^4}|Rm(g_i)|^{2}\leq \Lambda_1.\]

\remark In fact, under the controlled topology, i.e. $b_1(M)\leq
b_0$, the full Sobolev constant bound can be released to a lower
volume growth assumption, i.e. $\textrm{Vol}(B(x,r))\geq V_0r^n$.
Gang Tian and Jeff Viaclovsky's proof for Critical metric
\cite{TV08} also holds in our case, since the volume is continuous
in the $C^{1,\alpha}$ - topology.

Analogous convergence result also hold for the metric with bounded
Bach tensor and $C^2$ bound of the Scalar curvature, which is
actually easier to prove since we have the precise elliptic system
(\ref{BaSBC}).

\theorem \label{CTwBaSB} Let $\{M_i,g_i\}$ be a sequence of
Riemannian manifolds with
\begin{enumerate}
  \item either bounded Bach tensor, i.e., $|B(g_i)|\leq\Lambda$ or
  $\{M_i,g_i\}$ is  K\"{a}hler metric;
  \item $C^2$ bounded scalar
  curvature, i.e., $|\nabla^2S(g_i)|\leq\Lambda$;
  \item unit volume, i.e., $|\textrm{Vol}(M_i,g_i)|\equiv 1$;
  \item bounded Sobolev constant  $C_S$ in (\ref{USobIM});
\item bounded $L^{\frac{n}{2}}-$ norm of curvature, i.e.,
\[\{\int_M|Rm|^{\frac{n}{2}}\}^{\frac{n}{2}}\leq \Lambda_1.\]
\end{enumerate}
Then there exist a subsequence $\{j\}\subset\{i\}$  such that
$(M_j,g_j)$ converge to  a compact multi - fold
$(M_\infty,g_\infty)$ in the $C^{3,\alpha}$ topology off finite
singularity (as theorem \ref{MTCC}).

\example(\cite{Schoen87})There exists a family of conformally flat
metrics with constant scalar curvature $g_t$ on $S^1\times S^3$
($t\geq 1$) such that the diameter of $S^1\times\{x\}$ goes to $0$
as $t\rightarrow \infty$ for a point $x\in S^3$. For other point
$y\neq x$ the diameter of the slice $S^1\times\{y\}$ is bounded from
below. The limit space is the quotient space of $S^4$ which
identifies the north pole and the south pole.

\example (\cite{KT87} \cite{Nakajima92}) Let $M_\infty$ be an
orbifold given by the $\mathbb{Z}/2\mathbb{Z}$ - quotient of the
complex $2$- torus $T^2=\mathbb{C}^2/\mathbb{Z}^4$, where
$-1\in\mathbb{Z}/2\mathbb{Z}$ acts on $T$ by
\[(z_1,z_2)\;\textrm{mod}\;\mathbb{Z}^4\mapsto(-z_1,-z_2) \;\textrm{mod}\;\mathbb{Z}^4. \]
The flat metric on $T$ descends to an orbifold metric $g_\infty$ on
$M_\infty$, and it is an orbifold Ricci flat anti - self dual
metric. Moreover, $M_\infty$ has a complex manifold structure (with
singularities) since the $\mathbb{Z}/2\mathbb{Z}$ - action is
holomorphic. Let us take the minimal resolution $\pi: M \rightarrow
M_\infty$. The singularities are sixteen simple singularities of
type $A_1$. The minimal resolution $M$ is called a Kummer surface
and is an example of K3 surfaces. Let
$\mathcal{S}=\{x_1\cdots,x_{16}\}$ be the singular set, and let
$E_1,\cdots,E_{16}$ be the exceptional divisor. These are complex
submanifolds of $M$ biholomorphic to $\mathbb{C}P^1$ with the self -
intersection number $-2$. By the solution of the Calabi conjecture
we have a unique Calabi - Yau metric in each K\"{a}hler class, which
is automatically anti - self - dual. Take a K\"{a}hler class, and
there exist a sequence of  Ricci - flat Kahler metric $g_i$, as
follows: (1) the volume of $M$ with respect to $g_i$ is equal to
$1$; (2) the volume of the exceptional divisor $E_k$ is equal to
$\frac{1}{i}$ for $k = 1,\cdots,16$. It can be shown that the metric
$g_i$ converges to $\pi^*g_\infty$ over $M\setminus
\cup_{k=1}^{16}E_k$, but condition (2) forces the metric to become
degenerate along $E_k$ as $i\rightarrow\infty$, i.e. $E_k$ collapse
to a point, and the Riemannian curvature concentrates along
$\cup_{k=1}^{16}E_k$. Moreover, the curvature concentrates so
completely that the limit metric is a flat orbifold - metric.

\textbf{Acknowledgements}: The author would like to thank his
advisor Professor Gang Tian for suggesting this problem and constant
encouragement. The author would also like to thank Yalong Shi for
numerous suggestions which helped to improve the presentation.

\section{Self - Dual Weyl Curvature}
Riemannian geometry in dimension $4$ has some special features which
are not present in any other dimension. In dimensions $n\geq5$ the
group $SO(n)$ is simple, and the space of Weyl tensors is
irreducible under $SO(n)$, whereas in dimension $4$ we have
$SO(4)=SO(3)\times_{\mathbb{Z}_2}SO(3)$, and a decomposition of the
space of Weyl tensors into two $SO(4) -$ irreducible components. The
fact that $SO(4)$ is not simple is reflected at the Lie algebra
level in the the decomposition of the bundle of $2$ forms into self
- dual and anti - self - dual parts under the Hodge star operator.

Let $(M^4,g)$ be a oriented four dimensional Riemannian manifold.
The Hodge star operator $*$ associated to $g$ takes
$\Lambda:=\Lambda^2TM$ to itself and satisfying
$**=1$,\[\varphi\wedge\ast\psi=\langle\varphi,\psi\rangle dvol_g,\]
where $\langle \cdot, \cdot\rangle $ denotes the inner product on
$\Lambda$ induced by $g$. Then $\Lambda$ admits a decomposition of
the form
 \[\Lambda=\Lambda^{+}\oplus\Lambda^{-},\] where $\Lambda^{\pm}$ is the $\pm 1$
 eigenspace of $*$. The curvature operator $Rm$, viewed as an endomorphism on $\Lambda$, has the following matrix expression \cite{Besse87}:
\[Rm=\begin{bmatrix}
         W^++\frac{S}{12}\textrm{Id}_{\Lambda^+} & \overset{\circ}{Ric} \\
         (\overset{\circ}{Ric})^t& W^-+\frac{S}{12}\textrm{Id}_{\Lambda^-} \\
      \end{bmatrix}
\]
where $\overset{\circ}{Ric}\in \hbox{Hom}(\Lambda^+,\Lambda^-)$ is
the trace free part of the Ricci curvature, $S$ is the scalar
curvature, and $W^{\pm}\in S_0^2(\Lambda^\pm)$ is the (anti -) self
- dual part of the Weyl tensor, $S_0^2$ denote traceless symmetric
endomorphisms. If we denote the projection operator by
\[P_{\pm}:=\frac{1}{2}(1\pm *): \Lambda\rightarrow\Lambda^\pm,\]
then
\[W^{\pm}=P_{\pm}\circ Rm\circ P_{\pm}-\frac{S}{12}\textrm{Id}_{\Lambda^\pm}.\]

\definition Let $(M^4,g)$ be a compact Riemannian four dimensional Riemannian manifold, the metric $g$ is called anti - self - dual, if
$W^{+}=0(\Leftrightarrow *W=-W$).

We note that reversing the orientation transfers the self - dual
part to the anti - self - dual part. For anti - self - dual metric,
if we reverse the orientation, then $W^{-}(g)=0$ and $g$ is called
to be self dual ($W^{-}(\mathbb{C}P^2,J,\omega_{FS})=0$, i.e.
$(\overline{\mathbb{C}P^2},J,\omega_{FS})$ is anti - self - dual).

There are many interesting examples of anti - self - dual metrics.
First, the (anti -) self - duality of the metric is a conformally
invariant property. In particular, the conformal flat metrics will
be anti - self - dual. Second, a large number of anti - self - dual
metrics are K\"{a}hler metrics with zero scalar, since the self -
dual part of the Weyl tensor is given by \cite{Besse87}(Proposition
16.62)
\[W^{+}=\hbox{diag}(\frac{S}{6},-\frac{S}{12},-\frac{S}{12})
\]
In particular, a K3 surface with Calabi - Yau metric is anti - self
- dual. This also follows from Gauss - Bonnet formula and the
Hirzebruch signature formula, since with the canonical orientation
$\tau=-16$ and $\chi=24$.

\definition (Bach tensor) Let $(M^n,g)$ be a compact Riemannian manifold, the the Bach tensor is defined by
 \begin{equation*}
  B_{ij}= \frac{1}{n-3}\nabla^k\nabla^l W_{ikjl}+\frac{1}{n-2}R^{kl}W_{ikjl}
\end{equation*}
Such Bach tensor is symmetric, trace free and divergence free $2$
tensor. In dimension 4, Bach tensor arise as the Euler - Lagrange
equations of the functional on the $L^2-$ norm of the Weyl curvature
tensor. Using the Bianchi identities we may rewrite the Bach flat
(K\"{a}hler) metric  with constant scalar curvature  as an elliptic
system \cite{TV05a},
\begin{equation}\label{CRK}
    \left\{ \begin{aligned}
         \Delta Rm &= L(\nabla^2 Ric) + Rm\ast Rm \\
         \Delta Ric &=Rm \ast Ric
         \end{aligned} \right.
\end{equation}

In dimension 4, using Bianchi identity, Bach tensor can also be
written as $B_{ij}= 2\nabla^k\nabla^l W^+_{ikjl}+R^{kl}W^+_{ikjl}$,
thus anti - self - dual metric will be Bach flat. As we have
mentioned earlier, the K\"{a}hler metric on complex surface with
zero scalar curvature will be anti - self -dual, so it is also Bach
flat. For higher dimension, K\"{a}hler with constant scalar metric
will still satisfies the equation (\ref{CRK}).

Let $M$ be a closed oriented smooth manifold. A smooth Riemannian
metric $g$ on $M$ is a smooth section of the bundle $S^2T^*M$ of
positive definite symmetric $2 -$ tensors. The space $\mathcal{M}$
of all Riemannian metrics on $M$ is a convex open cone in
$\Gamma(S^2T^*M)$. A Riemannian metric is given locally by
functions, so we can define the (Sobolev) norm on $\mathcal{M}$ with
respect to some fixed metric, and the prescribed curvature condition
can be viewed as partial differential equation on $\mathcal{M}$.

With this viewpoint, we want to use a priori estimate of Elliptic
equation to study the convergence theory for the metric with bounded
self - dual Weyl tensor. Note that self - dual Weyl tensor $W^+(g)$
is equivariant under the action of the diffeomorphism transformation
and conformal change, but if we fix the gauge, the prescribe self -
dual Weyl tensor does form an elliptic system. In fact, it is now a
standard technique when we studying the geometry problems, such as
compactness of geometric Calculus of Variations(Yang - Mills
instanton, Einstein metric), DeTurck trick for existence theory of
geometric (Ricci) flow...... For anti - self dual structure, it is
well known that the local structure of  moduli spaces is controlled
by an elliptic deformation complex \cite{KK92} \cite{Tian08}. For
reader's convenience, we would like make it more clear in a PDE view
point (\ref{PWSGE}). As a consequence, we get a crucial estimates
for this paper, i.e. a priori $L^p$ estimate (\ref{LpEESPE}) for the
metric with bounded self - dual Weyl tensor and Scalar curvature.

\theorem \label{ESPEADS} Let $M$ be a closed oriented smooth $4 - $
manifold, we consider the following map (equation):
\begin{eqnarray}
 L: \mathcal{M}&\longrightarrow&S_0^2\Lambda^+\oplus C^\infty(M)\oplus TM\nonumber\\
  g&\longmapsto&(W^+(g),\; S(g), \;\tau_{g,\bar{g}}(id))\label{PWSGE}
\end{eqnarray}
where
$\tau_{\tau_{g,\bar{g}}}(id):=\textrm{tr}_g\nabla^{g\otimes\bar{g}}
d(id)$ is the tension field of the identity map $id:
(M,g)\longrightarrow (M,\bar{g})$, and $\bar{g}$ is some fixed
background metric.

Then the principle symbol $\sigma(L)$ of the linearized operator of
$L$ at $(x,g)$ is injective, which is given as follows: for any
$\xi\in T_xM$ and $h\in S_x^2TM$,
 \begin{equation}\label{SymbWSTM}
\sigma(L)(x,\xi)(h)=(h_{\Lambda^+},-\textrm{tr}(h)|\xi|^2+h(\xi,\xi),
h(\xi,\cdot)- \frac{1}{2}\textrm{tr}h\cdot\xi),
 \end{equation}
where  $h_{\Lambda^+}$ is defined by (\ref{SLWeyl}).

In particular, the equation (\ref{PWSGE}) on metric is an elliptic
system of partial differential equations of mixed order. And
consequently, the $L^p$ theory will holds:
\begin{equation*}
  \|g\|_{W^{2,p}}\leq C(\|g\|_{L^p}+\|W^+(g)\|_{L^p}+\|S(g)\|_{L^p}+\|\tau_{g,\bar{g}}\|_{W^{1,p}}).
\end{equation*}

\proof Recall the self - dual part of Weyl tensor  $W^{+}$ is
defined by \[W^{+}=P_{+}\circ Rm\circ
P_{+}-\frac{S}{12}\textrm{Id}_{\Lambda^+},\] then
\begin{eqnarray*}
  \delta W_g^{+}(h)&=&P_{+}\circ \delta Rm_g(h)\circ P_{+}+\delta P_{+g}(h)\circ Rm\circ P_{+}+P_{+}\circ Rm\circ \delta P_{+g}(h)\\
  &&-\frac{1}{12}\delta S_g(h)\textrm{Id}_{\Lambda^+}-\frac{S}{12}\delta
  \textrm{Id}_{\Lambda_g^+}(h)\\
&=&P_{+}\circ \delta Rm_g(h)\circ P_{+}-\frac{1}{12}\delta
S_g(h)\textrm{Id}_{\Lambda^+}+ \textrm{l.o.t}.
\end{eqnarray*}
On the other hand,  the first variation of curvature is
\cite{Besse87}:
\[
\begin{aligned}
  \delta Rm_g(h)&(X,Y,Z,U)=\frac{1}{2}[h(R(X,Y)Z,U)-h(R(X,Y)U,Z)]\\
    &+\frac{1}{2}[\nabla_{Y,Z}^2h(X,U)+\nabla_{X,U}^2h(Y,Z)-\nabla_{X,Z}^2h(Y,U)
-\nabla_{Y,U}^2h(X,Z)]
\end{aligned}
\]
Note that the scalar term $\delta S_g(h)\textrm{Id}_{\Lambda^+}$
does not contribute to the the principle symbol of the second
operator of $\delta W_g^{+}$, since $\delta W_g^{+}$ is traceless.
Therefore, by taking the traceless symmetric part, the principle
symbol of the second operator of $\delta W_g^{+}$ is given by
\begin{eqnarray*}
  \sigma(\delta W_g^{+})(x,\xi)(h)&=&\sigma(P_{+}\circ \delta Rm_g\circ P_{+})(x,\xi)(h)\\
  &&-\frac{1}{3}\textrm{tr}_{\Lambda^+}(\sigma(P_{+}\circ \delta Rm_g\circ P_{+})(x,\xi)(h))g_{\Lambda^+}\\
  &=&h_{\Lambda^+}
\end{eqnarray*}
where $h_{\Lambda^+}$ is defined as follows: for $e$, $e' \in
 \Lambda^+(x)$, and $\{e_i\}$ is an orthonormal basis of $\Lambda^+(x)$, we have
 \begin{equation}\label{SLWeyl}
h_{\Lambda^+}(e,e')=h(e(\xi),e'(\xi))-\frac{1}{3}\sum_{i=1}^3h(e_i(\xi),e_i(\xi))g_x(e,e').
 \end{equation}
The liberalization of the scalar curvature map is given by
\[\delta S_g(h)=\Delta \textrm{tr}h+\delta^2h-g(Ric, h),\]
then its symbol is\[\sigma(\delta
S_g)(x,\xi)(h)=-g(\xi,\xi)\textrm{tr}h+h(\xi,\xi).\] The
liberalization of the the tension field map is given by
\[\delta \tau_{g,\bar{g}}(id)(h)=\langle\nabla d(id),h\rangle-\langle\delta h+\frac{1}{2}d(\textrm{tr}h),d(id)\rangle,\]
therefore,
\[\sigma(\delta\tau_{g,\bar{g}(id)}g)(x,\xi)(h)=h(\xi,\cdot)-
\frac{1}{2}\textrm{tr}h\xi.\]

Combining the above three symbol computation gives (\ref{SymbWSTM}).

Now we will verify that the symbol of $\sigma(L)$ is injective, i.e.
$\forall \xi\neq 0$, $\sigma(L)(x,\xi)(h)=0$ implies that $h=0$. It
is easy to check that $\sigma(L)(x,\xi)(h)=0$ implies
$h_{\Lambda^+}=0$, $h(\xi,\xi)=0$, $\textrm{tr}h=0$.

First choose an orthonormal basis $X_1, X_2, X_3, X_4$ for $T_xM$,
then \[e_1=\frac{\sqrt{2}}{2}(X_1\wedge X_2+X_3\wedge
X_4),e_2=\frac{\sqrt{2}}{2}(X_1\wedge X_3+X_4\wedge
X_2),e_3=\frac{\sqrt{2}}{2}(X_1\wedge X_4+X_2\wedge X_3)\] gives an
orthonormal bases of $\Lambda_x^+$. For any $\xi=\xi^iX_i$ with
$|\xi|=1$, we have
\[
\begin{bmatrix}e_1(\xi) \\e_2(\xi)\\e_3(\xi)\\\frac{\sqrt{2}}{2}\xi\end{bmatrix}=
\frac{\sqrt{2}}{2}\begin{bmatrix}\xi^2&-\xi^1&\xi^4&-\xi^3 \\
\xi^3&-\xi^4&-\xi^1&\xi^2\\
\xi^4&\xi^3&-\xi^2&-\xi^1\\\xi^1&\xi^2&\xi^3&\xi^4\end{bmatrix}
\begin{bmatrix}X_1\\X_2\\X_3\\X_4\end{bmatrix}:=\frac{\sqrt{2}}{2}UX
\]
If $h_{\Lambda^+}=0$ and $h(\xi,\xi)=0$, from (\ref{SLWeyl}), then
it is equivalent to
\[
\begin{bmatrix}(h(e_i(\xi),e_j(\xi))&0\\0&\frac{1}{2}h(\xi,\xi)\end{bmatrix}=
\frac{1}{2}UhU^T=\frac{1}{3}\sum_ih(e_i(\xi),e_i(\xi)\begin{bmatrix}I&0\\0&0\end{bmatrix}
\]
Since $U\in O(4)$, then we have
\[h=\frac{2}{3}\sum_ih(e_i(\xi),e_i(\xi)U^T\begin{bmatrix}I&0\\0&0\end{bmatrix}U\]
Furthermore, if $\textrm{tr} h=0$, then
\[0=\textrm{tr} h=2\sum_ih(e_i(\xi),e_i(\xi),\]
and consequently, $h=0$.

The symbol $\sigma(L)$ is injective,  by comparing the dimension
(\textrm{dim} =10), we conclude that the equation (\ref{PWSGE}) is
an elliptic system of partial differential equations with mixed
order. For an elliptic system (of mixed order), we have a priori
estimate, for example, the Schauder or $L^p$ theory \cite{Morrey66}.
Alternatively, if we replace the tension field
$\tau_{g,\bar{g}}(id)$ in equation (\ref{PWSGE}) by $\mathcal
{L}_{\tau_{g,\bar{g}}(id)}g$, then it will be a overdetermined
elliptic system of second order partial differential equations. In
fact, the a priori estimate ($L^p$ theory)  hold for differential
operators between vector bundles is equivalent to the injectivity of
the symbol, while the `solubility criteria' holds if the symbol is
surjective (theorem 19.25 in \cite{Palais68} or theorem A.8 in
\cite{DK90}).

In local harmonic coordinates, since the metric is $C^{1,\alpha}$
close to the Euclidean metric, if necessary, we can make the
harmonic norm $C$ small enough, then the operator $L$ is uniformly
elliptic with $C^{1,\alpha}$ continuous leading terms, by the $L^p$
theory, we have the estimate
\begin{eqnarray}
  \|g\|_{W^{2,p}(B_r)}&\leq&C(\|g\|_{L^p(B_{2r})}+\|W^+(g),S(g)\|_{L^p(B_{2r})}+\|\tau_{g,\bar{g}}(id)\|_{W^{1,p}(B_{2r})})\nonumber\\
  &=&C(\|g\|_{L^p(B_{2r})}+\|W^+(g)\|_{L^p(B_{2r})}+\|S(g)\|_{L^p(B_{2r})})\label{LpEESPE}
\end{eqnarray}
The last step holds, since $\tau_{\tau_{g,\bar{g}}}(id)=0$ if we
identify the geodesic ball with the Euclidean ball under harmonic
coordinate. \qed

\section{Convergence Theory}

To generalize the convergence theory with assumption on sectional
curvature to weaker curvature bounded hypothesis, it is convenient
to use the concept of harmonic radius, which was introduced and
developed in \cite{Anderson90} \cite{AnCh92} \cite{Petersen97}.
\definition (Harmonic Radius)  Let $(M,g)$ be a Riemannian manifold, fixing any $m\in\mathbb{N}$,
given any $x\in M$, there is $r=r(x)$, for which we can find a harmonic coordinate system
\[(\{x_i\}_i^n) : B(x, r)\subset (M,g)\rightarrow\mathbb{R}^n\]
such that
\begin{enumerate}
  \item the coordinate function $(\{x_i\}_i^n):B(x, r)\subset (M,g)\rightarrow\mathbb{R}^n$ is harmonic;
  \item the metric tensor $g_{ij}:=g(\nabla x_i,\nabla x_j)$ is $C^{m,\alpha}$ bounded on $B(x, r)$, i.e.
   \[e^{-C} \delta_{ij}\leq g_{ij}\leq e^C \delta_{ij}\;(\textrm{as bilinear forms}),\]
   and
  \[\sum_{1\leq|\beta|\leq m} r^{|\beta|}\sup|\partial^\beta g_{ij}|+\sum_{|\beta|= m}r^{m+\alpha}[\partial^\beta g_{ij}]_{\alpha}\leq C,\]
\end{enumerate}
for some fixed constant $C\geq 0$, where the norms are taken with
respect to the coordinates $(\{x_i\}_i^n)$ on $B(x, r)$. We say that
$x\in (M,g)$ admits a harmonic coordinate with bounded $C^{m,\alpha}
- $ norm on the scale $r$:
\[\|x\in(M,g)\|_{C^{m,\alpha},r}\leq C.\] Moreover, we let $r_h(x)$ be the $C^{m,\alpha}$ harmonic radius at
$x$, which is defined as the radius of the largest geodesic ball
about $x$, on which there are $C^{m,\alpha}$ harmonic coordinates,
i.e.
\[r_h(x)=\sup\{r>0|\;\|x\in(M,g)\|_{C^{m,\alpha},r}\leq C\}.\]
The $C^{m,\alpha}$ harmonic radius of $(M,g)$ is defined by
$r_h(M)=\inf_{x\in M}r_h(x)$.

\definition \label{MCCN} For a compact metric space $X$, define the covering number of the geodesic ball on the scale $\epsilon$ as follows
\[\textrm{Cov}(\epsilon)= \min\{n|\exists\{x_i\}_{i=1}^n\subset X, \bigcup_{i=1}^nB(x_i, \epsilon)=X,  B(x_i,\frac{\epsilon}{2})\bigcap_{i\neq j} B(x_j,\frac{\epsilon}{2})=\emptyset\}\]
Now,  let us state the fundamental theorem of convergence theory:
\theorem \label{FTCT} (\cite{Anderson90} \cite{Petersen97}) For
given
 $n\geq 2$, $C\geq 0$, $N>0$, $\alpha \in (0,1]$ and $r_0>0$, consider the class
$\mathcal{M}(n,C,N,r_0)$ of $n - $ dimensional Riemannian manifolds
\[\{(M,g)|r_h(x)\geq r_0, \|x\in(M,g)\|_{C^{m,\alpha},r_0}\leq C, \forall x \; \hbox{and}\; \textrm{Cov}(\frac{r_0}{10})\leq N.\}\]
Then $\mathcal{M}(n,C,N,r_0)$ is compact in the $C^{m, \beta}$
Cheeger - Gromov topology for all $\beta<\alpha$. Moreover, the
theorem also valid for bounded domains in Riemannian manifolds, as
well as for pointed complete Riemannian manifolds, provided one
works on compact subsets.

\theorem \label{ERTm}(\cite{TV05a} \cite{TV08} \cite{Carron10}) Let
$(M,g)$ be a complete Riemannian manifold or Riemannian multi - fold
with finite point singularities, $g$ is critical metric (for
example, Bach flat with zero scalar curvature). Assume that
$(M^n,g)$ satisfies the Sobolev inequality,
\[\{\int_M|v|^{\frac{2n}{n-2}}d\nu\}^{\frac{n-2}{n}}\leq C_S\int_M|d v|^2d\nu, \forall v\in C_0^{0,1}(M);\]
and the curvature has bounded $L^{\frac{n}{2}}$ norm,
\[\{\int_M|Rm|^{\frac{n}{2}}d\nu\}^{\frac{n}{2}}\leq \epsilon.\]
Then $(M^n,g)$ has finite many ends, which is ALE of order two.
Moreover, the volume is at most Euclidean volume growth,
$\textrm{Vol}(B(p,r))\leq Vr^n$ for some positive constant
$V=V(n,C_S, \epsilon)$.

Moreover, if $(M,g)$ is smooth and $\epsilon=\epsilon_0$ is small
enough, which depends on $n$ and the Sobolev constant $C_S$, then
$(M,g)$ is isometric to the Euclidean space $(\mathbb{R}^n,g_E)$.

\lemma \label{HaRiEB} For any $m\geq 2$, let $B(r):=B(x_0,r)$ be a
geodesic ball in a compact oriented four Riemannian manifold
$(M^4,g)$, where $g$ has bounded self - dual Weyl tensor (or
K\"{a}hler) and bounded scalar curvature, $|W^+|+|S|\leq \Lambda$.
Then there exist a positive constant $\epsilon_0=\epsilon(C_S)$ and
$\kappa_0=\kappa_0(C_S,\Lambda)$ such that if
\begin{equation}\label{EsEB}
\{\int_{B(x,2r)}|Rm|^{2}d\nu\}^{\frac{1}{2}}\leq \epsilon_0,
\end{equation}
then for all $x \in B(x_0,r)$, the $C^{1, \alpha}$ harmonic radius
$r_h(x)$ satisfies
\begin{equation}\label{LHRE}
\frac{r_h(x)}{\hbox{dist}(x,\partial B)}\geq \kappa_0>0.
\end{equation}

\proof On a smooth fixed smooth Riemannian manifold $(M,g)$, it is
clear that the harmonic radius $r_h(x)$ is positive or saying
(\ref{LHRE}) holds, but $\kappa_0$ depends on $(M,g)$ and $x$. Thus,
we must show that $\kappa_0$ depends only on the hypothesis
prescribed in the lemma.

We argue by contradiction, which is similar as the blow up analysis
for Ricci curvature case in \cite{Anderson90}. If (\ref{LHRE}) is
false, then there are sequence of Riemannian $4$ manifolds
$\{(M_i,g_i)\}$ with the bounds in the lemma, and points $x_i\in
B_i(r)\subset (M_i,g_i)$ such that
\begin{equation}\label{HaRaRa}
\frac{r_h(x_i)}{\hbox{dist}(x_i,\partial B_i)}\rightarrow 0,\quad
\hbox{as}\; i\rightarrow \infty.
\end{equation}
We may assume, with out lose of generality, that the points $x_i$
realize the minimum of the left side of (\ref{LHRE}) and
\[\|x_i\in(M_i,g_i)\|_{C^{1,\alpha},r_h(x_i)} \in [\frac{C}{2}, C].\]
By scaling theses metrics suitably, namely,
$\bar{g}_i=r_h(x_i)^{-2}g_i$, then
\begin{enumerate}
  \item $\bar{r}_h(x_i)=1$ and $\bar{r}_h(x)$ is bounded below on balls of finite distance to $x_i$,
  which follows from scale invariant property of Harmonic
  norm \cite{Petersen97},
\[\|x\in(M,\lambda^2g)\|_{C^{m,\alpha},\lambda r}=\|x\in(M,g)\|_{C^{m,\alpha},r};\]
  \item $\hbox{dist}_{\bar{g}_i}(x_i,\partial B_i)\rightarrow \infty$, since the ratio in (\ref{HaRaRa}) is scale invariant;
  \item $|W^+(\bar{g}_i)|+|S(\bar{g}_i)|\leq r_h^2(x_i)\Lambda\rightarrow 0$, and the curvature have $\epsilon- $ small $L^{2}$  norm \[C_S\{\int_{B(x,\frac{2r}{r_h(x_i)})}|Rm(\bar{g}_i)|^2d\bar{\nu}_i\}^{\frac{1}{2}}\leq \epsilon_0,\]
\end{enumerate}
with respect to the metric $\bar{g}_i$.

We now consider, the sequence of pointed Riemannian manifolds
\[\{(B_i(x_i, \frac{r}{r_h(x_i)}), x_i, \bar{g}_i)\subset
(M_i,x_i,\bar{g}_i)\},\] by the fundamental theorem of convergence
theory \ref{FTCT}, then the sequence subconvergent, in the pointed
$C^{1,\beta}$ topology $(\forall \beta<\alpha)$, uniformly on
compact subsets, to a complete $C^{1,\alpha}$ Riemannian manifold
$(N, \bar{x}, h)$.

\claim The convergence is actually better, namely in the
$C^{1,\alpha}$ topology, where $\alpha$ is given by the hypothesis
of the lemma.

Moreover, we can even prove more than we need, i.e. the convergence
is in the $C^{1,\alpha}\cap W^{2,p}$ topology, for any $\alpha<1$
and $1<p<\infty$. By Soblev embedding theorem, $W^{2,p}\subset
C^{1,\alpha}$ if ${p>n}$, so it suffices to prove the convergence is
in the $W^{2,p}$ topology. To see this, by theorem \ref{ESPEADS}, we
know that the prescribed self - dual Weyl tensor and scalar
curvature equation (\ref{PWSGE}) is an elliptic system of partial
differential equations of second order under harmonic coordinate:
\begin{equation} \label{ASDEPSC}
\left\{ \begin{aligned}
         W^{+} &=L(g^{-1}\partial\partial g)+Q_1(\partial g,\partial g)\in L^\infty\\
          S (g) &=-\frac{1}{2}g^{ij}g^{kl}\frac{\partial^2}{\partial x^k\partial x^l}g_{ij}+Q_2(\partial g,\partial
          g)\in L^\infty
                          \end{aligned} \right.
                          \end{equation}
where $L$ denotes linear combination,  $Q$ is a quadratic term in
the first order derivatives of $g$. For a priori estimate, since
$\|g_{ij}-\delta_{ij}\|_{C^{1,\alpha}}<C$, if necessary we can make
$C$ small, the above system actually can be viewed as a uniform
linear elliptic system of $g_{ij}$ with $C^{1,\alpha}$ coefficients.
By the a priori estimate (\ref{LpEESPE}), the $L^p$ theory for
elliptic systems gives a uniform bound on $\|g\|_{W^{2,p}}$ for any
$1<p<\infty$,
\[\|g\|_{W^{2,p}}\leq C(\|g\|_{L^{p}}+\|Q(\partial g,\partial g)\|_{L^{p}}+\|W^+(g)\|_{L^p}+|S(g)\|_{L^p})\leq C.\]
As a consequence, the convergence is in the $C^{1,\alpha}\cap
W^{2,p}$ topology, for any $\alpha<1$ and $1<p<\infty$.

More precisely, the $(B_i(\frac{r}{r_h(x_i)}),x_i,\bar{g}_i)$ are
covered by harmonic coordinates that converge in the $C^{2,\alpha}$
topology to the harmonic coordinates on limit space $N$, and the
metric coefficients $\bar{g}_i$ converge in the $C^{1,\alpha}$
topology to $h$.

Since the $C^{1,\alpha}$ norm is continuous and harmonic radius is
continuous with respect to $C^{1,\alpha}$ or $W^{2,p}$ convergence
\cite{Anderson90} \cite{Petersen97}, and consequently,
\begin{equation}\label{HRx1}
r_h(\bar{x})= 1,\;
\|\bar{x}\in(N,h)\|_{C^{1,\alpha},r_h(\bar{x})}\geq \frac{C}{2}>0.
\end{equation}

\claim $h$ is a smooth Riemannian  metric and $(N,h)$ is isometric
to the Euclidean space $(\mathbb{R}^4,g_E)$.

Since the convergence is in the $C^{1,\alpha}\cap W^{2,p}$ topology,
we can therefore conclude that the limit metric $h$ is a weak
$C^{1,\alpha}\cap W^{2,p}$ solution of the elliptic system,
\begin{equation}
\left\{ \begin{aligned}
         W^{+} &=L(g^{-1}\partial\partial g)+Q_1(\partial g,\partial g)=0\\
          S (g) &=-\frac{1}{2}g^{ij}g^{kl}\frac{\partial^2}{\partial x^k\partial x^l}g_{ij}+Q_2(\partial g,\partial
          g)=0
                          \end{aligned} \right.
                          \end{equation}
namely, the anti self dual or K\"{a}hler metric with zero scalar
curvature is a second order divergence form qusi-linear elliptic
system of the metric modulo diffeomorphisms by theorem
\ref{ESPEADS}. With the a priori estimate (\ref{LpEESPE}) and a
standard bootstrap argument, and also Sobolev embedding theorem, we
conclude that the metric $h$ is actually a smooth (in fact,
analytic) Riemannian metric with
\begin{equation*}
C_S\{\int_{N}|Rm(h)|^{2}d\nu_h\}^{\frac{1}{2}}\leq \epsilon_0.
\end{equation*}

If $\epsilon=\epsilon_0$ is sufficiently small (which will depends
only on the Sobolev constant), by the $\epsilon - $ rigidity theorem
\ref{ERTm}, we conclude that $Rm(h)\equiv 0$, i.e. $N$ is flat. On
the other hand, bounded Sobolev constant implies Euclidean volume
growth. Consequently, $(N,h)$ is isometric to the Euclidean space
$(\mathbb{R}^n,g_E)$.

Since the Euclidean space admits global harmonic coordinates, i.e.
\[r_h(x)=\infty, \quad \|x\in(\mathbb{R}^n,g_E)\|_{C^{1,\alpha},r}=0,\; \forall\; r>0.\]
However, this violates (\ref{HRx1}). \qed

Now we can prove the main theorem \ref{MTCC}, which  is an immediate
consequence of the main lemma \ref{HaRiEB} on harmonic radius
estimate and the fundamental theorem \ref{FTCT} of convergence
theory.

\theorem \label{CoLiTh} With the same hypothesis in theorem
\ref{MTCC}, then there exist a subsequence $\{j\}\subset\{i\}$  such
that $(M,g_j)$ converge to  a compact metric space
$(M_\infty,g_\infty)$ in the $C^{1,\alpha}$ topology outside the
finite singular set $S=\{x_1,\cdots, x_s\}$.

\proof As in the case of bounded Ricci curvature or or Bach flat
metric with constant scalar curvature, take $\epsilon=\epsilon_0$ in
theorem \ref{ERTm}, consider the sets
\[\mathcal {R}_i(r)=\{x\in M_i|\{\int_{B(x,2r)}|Rm|^{\frac{n}{2}}\}^{\frac{n}{2}}<\epsilon_0\}\]
and
\[\mathcal {S}_i(r)=\{x\in M_i|\{\int_{B(x,2r)}|Rm|^{\frac{n}{2}}\}^{\frac{n}{2}}\geq \epsilon_0\}\]
then $M_i=\mathcal {R}_i(r)\cup\mathcal {S}_i(r)$, and also
$\mathcal {R}_i(r_1)\subset\mathcal {R}_i(r_2)$, $\mathcal
{S}_i(r_1)\supset\mathcal {S}_i(r_2)$, for any $r_1 > r_2$.

For all $x\in\mathcal {R}_i(r)$, by the main lemma \ref{HaRiEB}, we
have the estimate on  $C^{1,\alpha}$ harmonic radius,
\[r_h(x)\geq \kappa_0r,\] where $\kappa_0=C(C_S,\Lambda)$.
On the other hand, the uniform Sobolev constant implies
noncollapsing, namely, $\textrm{Vol}(B(x,r))\geq C(C_S)r^n$. Then
the covering number (see definition \ref{MCCN}) on any compact
subset of $\mathcal {R}_i(r)$ on the scale $\inf_{x\in\mathcal
{R}_i(r)} r_h(x)\geq \kappa_0r$ can be bounded by
\[\textrm{Cov}(\frac{1}{10}\kappa_0r)\leq\frac{\textrm{Vol}(M)}{\textrm{Vol}(B(x,\frac{1}{10}\kappa_0r))}\leq
\frac{C(C_S,\Lambda)}{r^4}.\]

 With the fundamental convergence theorem \ref{FTCT}, the sequence $(\mathcal {R}_i(r),g_i)$ is $C^{1,\alpha}$
subconvergent
 to a $C^{1,\alpha}$ (open) Riemannian
manifold $(\mathcal {R}_\infty(r),g_\infty)$ on the compact set.

To construct the limit space,  we will be brief since this step is
quite standard, see for example \cite{Anderson89} \cite{BKN89}
\cite{Tian90},  and also \cite{Anderson05} \cite{TV05b}.

We now choosing a sequence $\{r_j\}\rightarrow 0$ with
$r_{j+1}<\frac{1}{2}r_j$, pepeat the above construction by choosing
subsequence, we still denote $\{j\}$. Since $\mathcal
{R}_i(r_{j})\subset\mathcal {R}_i(r_{j+1})$, then we have a sequence
of limit spaces with natural inclusions
\[\mathcal{R}_\infty(r_{j})\subset\mathcal {R}_\infty(r_{j+1})\subset\cdots\subset\mathcal {R}_\infty:=\textrm{dir.}\lim\mathcal{R}_\infty(r_{j})\]
By the $C^{1,\alpha}$ convergence, $(\mathcal {R}_\infty,g_\infty)$
is $C^{1,\alpha}$ (open) Riemannian manifold, and there are
$C^{2,\alpha}$ smooth embedding $F_i: (\mathcal
{R}_\infty,g_\infty)\rightarrow (M,g_i)$ such that
$F_i^*g_i\rightarrow g_\infty$ in the $C^{1,\alpha}$ topology on any
compact set of $\mathcal {R}_\infty$.

Letting $\{B(x_k^i, \frac{r}{4})\}_{k\in\mathbb{N}},
r<\frac{1}{4}\rho_0$, where $\rho_0$ is the Euclidean volume growth
scale in theorem \ref{UEUVSE},  be a collection of a maximal family
of disjoint geodesic balls in $M_i$, then $M_i\subset \cup_kB(x_k^i,
r)$. There is a  uniform bound, independent of $i$,  on the number
of points $\{x_k^i\in {S}_i(r)\}$, which follows from
\begin{equation}\label{FpSS}
m\leq\sum_{i=1}^m\epsilon_0^{-2}\int_{B(x_k^i,2r)}|Rm|^{2}\leq
C\epsilon_0^{-2}\int_{M_i}|Rm|^{2},
\end{equation}
where $C=\sup_{x\in
M_i}\frac{Vol(B(x,\frac{9r}{4})}{Vol(B(x,\frac{r}{4})}\leq
C(C_S,\Lambda)$. The last inequality holds because we have upper
bound of the volume growth (\ref{UVGE}).

Let $\{r_j\}$ be as above without loss generality, we will assume
$m$ is fixed, i.e. the number of mutually disjoint balls, which is
centered in $S_i(r_j)$ and has radius $\frac{r_j}{4}$, is
independent on $i$ and $j$. As a consequence, every point of
$S_i(r_j)$ is contained in a ball of diameter no greater than
$mr_j$. Hence, most of the volume $(M_i,g_i)$ is contained in
$\mathcal{R}(r_j)$.  Using the embedding $F_i^j: (\mathcal
{R}_\infty(r_j),g_\infty)\rightarrow(\mathcal {R}_i(r_j),g_i)$, we
see that for any fixed $j$, and $i$ sufficiently large, arbitrarily
large compact subsets od ${R}_\infty\setminus {R}_\infty(r_j)$ are
almost isometrically embedded into $m$ disjoint balls of radius
$r_j$. Letting $j\rightarrow \infty$, it follows that the the
boundary components shrink to points with  respect $g_\infty$. In
other words, one can add finite points $\mathcal
{S}_\infty=\{x_1,\cdots,x_m\}$ to $\mathcal {R}_\infty$ such that
$M_\infty:=\mathcal {R}_\infty\cup\mathcal {S}_\infty$ is complete
with respect to the length structure $g_\infty$, i.e. the Riemannian
metric has a $C^0$ extension across the singularity.  Moreover,
$(M_i,g_i)$ sub convergent to $M_\infty$ in the Gromov - Hausdorff
topology and the volume of geodesic ball (may contain singularity)
is continuous with respect the $C^{1,\alpha}$ convergence (off
finite singularity).

We now examine the topological structure near the singularity by
essential studying the tangent cones at the singularity. Fix $p\in
S_\infty\subset M_\infty$, let $r(x)=\textrm{dist}(x,p)$ and denote
the annulus around $p$ to be
\[A(r_1,r_2)=\{x\in M_\infty|r_1<r(x)<r_2\},\] where
$r_1<r_2<\textrm{dist}(p,S_\infty\setminus\{p\})$. By the
$C^{1,\alpha}\cap W^{2,p}$ convergence, recall the $C^{2,\alpha}$
smooth embedding $F_i: (\mathcal {R}_\infty,g_\infty)\rightarrow
(M,g_i)$, the curvature will converge in the $L^p -$ sense, and then
\begin{equation*}\int_{F_i(\mathcal
{R}_\infty)}|Rm(g_i)|^2<\infty,\;\forall i.
\end{equation*}
 In particular, for $\epsilon_0$ in theorem \ref{ERTm}, there is an $r_0>0$ such that
\begin{equation}\label{CNATZ}
\int_{F_i(A(0,r_0))}|Rm(g_i)|^2\leq  \epsilon_0^2, \forall i.
\end{equation}

Now we do blow up analysis on $M_\infty$, it is equivalent to blow
up the sequence. Namely, given any sequence $r_j\rightarrow 0$, let
$j\rightarrow 0$, the metric annulus $(F_i(A(\frac{s_j}{j},
js_j)),\frac{1}{s_j^2}g_i)$ (by taking diagonal sequence) sub -
converge to $C^{1,\alpha}\cap W^{2,p}$ annulus $(A_\infty(0,
\infty), g_\infty)$, where $g_\infty$ weak solution of anti self -
dual with zero scalar curvature equation. With the regularity theory
of Elliptic equation, it follows that $(A_\infty(0, \infty),
g_\infty)$ is smooth. On the other hand, by(\ref{CNATZ}), we know
\begin{equation}
\int_{A(0,\infty)}|Rm(g_\infty)|^2\leq  \epsilon_0^2,
\end{equation}
with Sobolev constant, we conclude that each component of
$(A_\infty(0, \infty), g_\infty)$ is isometric to the Euclidean cone
on a space form $S^3/\Gamma$ for some finite subgroup of $O(4)$.

If one has lower Ricci curvature, then the limit orbifold is
irreducible, which is proved in \cite{Anderson89} by means of the
Cheeger - Gromoll splitting theorem. In our case, there may be more
than one cones associated to one singularity. If we perform a
standard bubble analysis, one can estimate the precise bound on the
end of associated ALE space, which in turn implies a bound on the
number of cones at each singular point, depending only on
$\|Ric_-\|_{L^2}$, $C_S$, see \cite{Carron98} and \cite{TV08}. This
also give a alternative way to show the limit orbifold is
irreducible if one has lower Ricci curvature. For K\"{a}hler metric,
only irreducible singular points can occur in limit, i.e. orbifold
point, see more details in \cite{TV05b}.

It follows that the neighborhoods of each singular points is
homomorphic to finite cones on spherical spaces forms. \qed

\remark  For the proof of theorem \ref{CTwBaSB}, the argument is
similar. In fact, it is much easier to estimate the harmonic radius
as did in the main lemma \ref{HaRiEB}.

With Bianchi identity, the laplacian of Ricci curvature is related
to the Bach tensor (K\"{a}hler) and Scalar curvature \cite{TV05a},
so we have a coupled system:
 \begin{equation}\label{BaSBC}
 \left\{ \begin{aligned}
                  \Delta Ric &=2B+\frac{1}{3}\hbox{Hess} S+Rm*Ric\\
                  \Delta g&=Q(\partial g, \partial g)-2Ric
                           \end{aligned} \right.
                           \end{equation}
Under the $C^{3,\alpha}$ harmonic coordinates, we have the improved
estimate \[\|g\|_{C^{3,\alpha'}}<C,\;\forall
0<\alpha'<1;\;\textrm{and}\;\|g\|_{W^{4,p}}<C,\;\forall
1<p<\infty.\] Moreover, the blow up limit will be flat since
$\epsilon - $ rigidity theorem \ref{ERTm} hold for Bach flat
(K\"{a}hler) metric with zero Scalar curvature.  The left argument
is similar and will be omitted here.

\section{Volume Growth Near Singularity}
We have already seen that the volume growth plays a crucial role in
understanding the structure near the singular set, see (\ref{FpSS}).
By lack of the volume comparison, we must find alternative approach
to bound the volume on a fixed scale, i.e. for some $\rho>0$, there
exist $V_1>0$ such that $ \textrm{Vol}(B(p,r))\leq V_1r^n, \forall
r<\rho$. For Bach flat metric with constant Scalar curvature, Gang
Tian and Jeff Viaclovsky concluded that the volume does bound on all
scale,  and the bound depends  only on the Sobolev constant and
$L^{2}$ norm of curvature \cite{TV05a} \cite{TV05b} \cite{TV08}. In
fact, if we check their paper carefully, we will find that the
argument also holds in our case, where we work in $C^{1,\alpha}$
category in place of the $C^\infty$ category. The difficulty caused
by the concentrate of curvature. If we do blow up analysis
carefully, as did in Einstein case \cite{Nakajima92}, there will
bubble out some non - flat ALE space (tree) which will satisfy
stronger geometric conditions; and consequently, we can bound the
volume growth.

 \theorem\label{UEUVSE} Let $(M^4,g)$ be a
compact oriented four manifold with bounded self - dual Weyl tensor
and Scalar curvature, i.e. $|W^{+}(g)|+|S(g)|\leq \Lambda$; bounded
Sobolev constant $C_S$; and also finite $L^2$ - curvature, i.e.
$\|Rm\|_{L^2}<\Lambda_1$. For some $\rho_0>0$, there exists a
constant $V_1>\omega_4$, depending only upon $\Lambda$, $\Lambda_1$,
$C_S$ such that
\begin{equation}\label{UVGE}
 \textrm{Vol}(B(x,r))\leq V_1r^4
\end{equation}
for all $x\in M$ and $0<r<\rho_0$.

\proof The theorem can be established by proceeding the same bubble
procedure in \cite{TV08}. For reader's convenience, we will copy
down their argument with some slight modification to give a full
argument in our case.

In the first place, if the curvature does concentrate too much, i.e.
for $\rho>0$,
\[\int_{B(x,2\rho)}|Rm|^2<\epsilon_0^2,\]
the volume growth will be controlled. In fact, by lemma
\ref{HaRiEB}, the harmonic radius of $B(x,\rho)$ is bounded below,
namely, there is uniform constant $r_0$, such that, for all $y\in
B(x,\rho)$,
\[r_h(y)>r_0\rho,\quad \|y\in(M,g)\|_{C^{1,\alpha},r_0\rho}<C,\] and
consequently,
\begin{equation}\label{UVGSE}
\textrm{Vol}(B(x,r))\leq e^{2C}\omega_4r^4,\; \forall r\leq \rho.
\end{equation}

For any metric $(M,g)$, define the maximal volume ratio on the scale
$\rho$ as
\[\textrm{MV}(g,\rho)=\max_{x\in M, 0<r<\rho}\frac{\textrm{Vol}(B(x,r))}{r^4}.\]
Note that for any compact smooth four Riemannian manifold $(M,g)$,
\[\lim_{\rho\rightarrow 0}\textrm{MV}(g,\rho)=\omega_4,\]
where $\omega_4$ is the volume ratio of the Euclidean metric on
$\mathbb{R}^4$.

In this paper, we consider the maximal volume ratio on finite scale
rather than on all scale in \cite{TV08}. On the one hand, the local
non - inflated volume is enough to shrink the singular set to point;
On the other hand, one will see, with lacking of $\epsilon$ -
regularity, we can not prove the Euclidean volume growth on the
large scale by volume comparison.

If the theorem is not true, then for any a sequence
$\rho_j\rightarrow0$ with $\rho_{j+1}<\frac{1}{2}\rho_j$, if we fix
$j$, there exists a sequence of metrics $(M,g_{i,j})$, which
satisfies the hypothesis in the theorem, but $MV(g_{i,j},
\rho_j)\rightarrow\infty$. By passing to a diagonal subsequence, for
any a sequence $\rho_i\rightarrow0$ with
$\rho_{i+1}<\frac{1}{2}\rho_i$,  there exists a sequence of metrics
$(M,g_{i})$, which satisfies the hypothesis in the theorem, but
\begin{equation}\label{CUVAC}
\textrm{MV}(g_{i}, \rho_i)\rightarrow\infty,\;\textrm{ as}\;
i\rightarrow\infty.
\end{equation}

For this sequence, we can exact a subsequence (which for simplicity
we continue to denote by the index $i$) and $r_i<\rho_i$ such that
\begin{equation}\label{VoRaBU}
2e^{2C}=\frac{\textrm{Vol}(B(x_i,r_i))}{r_i^4}= \max_{r\leq
r_i}\frac{\textrm{Vol}(B(x_i,r_i))}{r^4},
\end{equation}
where $e^{2C}$ comes from (\ref{UVGSE}).

We furthermore assume $x_i$ is choose so that $r_i$ is minimal, that
is, the smallest radius such that
\[\textrm{Vol}(B_{g_i}(x,r))\leq 2e^{2C}r^4, \forall\;x\in M_i\; \textrm{and}\; r\leq r_i.\]
Note that the inequality
\begin{equation}\label{CuCoSB}
 \int_{B(x_i,2r_i)}|Rm(g_i)|^2 \geq \epsilon_0^2
\end{equation}
must hold, each ball with lager volume growth  (singularity) takes
at least $\epsilon_0$ of $L^2$ - curvature. Otherwise, by the
estimate (\ref{UVGSE}), we would have
\[\textrm{Vol}(B_{g_i}(x_i,r_i))\leq e^{2C}r_i^4,\] which violates the
choice of $r_i$ in (\ref{VoRaBU}).

Now, we consider the resealed metric $\tilde{g}_i=r_i^{-2}g_i,$ so
that $B_{g_i}(x_i,r_i)=B_{\tilde{g}_i}(x_i,1)$. From the choice of
$x_i$ and $r_i$, the rescaled metrics $\tilde{g}_i$ have bounded
volume ratio, in all of unit size.

From the main theorem  \ref{MTCC}, their exist a subsequence
converges on compact subsets to a complete length space $(M_\infty,
g_\infty, x_\infty)$ in the $C^{1,\alpha}$ topology off finite many
singularities, where $(M_\infty, g_\infty, x_\infty)$ is a
multi-fold, $g_\infty$ is a smooth anti - self - dual metric with
zero Scalar curvature. Further, from theorem \ref{ERTm} for the
multi- fold case, see proposition 4.3 and claim 4.4, p.14 in
\cite{TV08}, there exists a constant $A_1$ such that
\begin{equation}\label{UVGLMD}
\textrm{Vol}(B_{g_\infty}(x_\infty,r))\leq A_1r^4, \;\textrm{for
all}\; r>0.
\end{equation}
One have seen that if $r_i\rightarrow 0$, then the blow up limit
will be a smooth multi - fold with critical metric, which is crucial
to conclude (\ref{UVGLMD}). This is the main reason that we consider
the maximal volume ratio on the finite scale.

We next return the (sub)sequence $(M,g_i)$ and exact another
subsequence so that
\begin{equation}\label{AVoRaA1}
2600 A_1=\frac{\textrm{Vol}(B(x_i',r_i))}{r_i'^{4}}=\max_{r\leq
r_i'}\frac{\textrm{Vol}(B(x_i',r))}{r^{4}}.
\end{equation}
Again, we assume that $x_i'$ is chosen so that $r_i'$ is minimal,
that is, the smallest radius for which
\[\textrm{Vol}(B_{g_i}(x,r))\leq 2600A_1r^4, \forall\;x\in M_i\; \textrm{and}\; r\leq r_i'.\]
Clearly, $r_i<r_i'<\rho_i\rightarrow0$.

Arguing as above, we repeat the rescaled limit construction, but now
with scaled metric $g_i'={r_i'}^{-2}g_i$, and basepoint $x_i'$. We
find a limiting multi - fold $(M_\infty',g_\infty',x_\infty')$, and
a constant $A_2\geq 2600 A_1$ so that
\[\textrm{Vol}(B_{g_\infty'}(x_\infty',r))\leq A_2r^4, \;\textrm{for all}\; r>0.\]
For the same reason as in (\ref{CuCoSB}), we must have
\[\int_{B_{g_i}(x_i,2r_i')}|Rm(g_i)|^2\geq\epsilon_0^2.\]

Since the $L^2$ - curvature is finite, and each lager volume growth
ball (singularity) takes at least $\epsilon_0$ of $L^2$ - curvature,
it is reasonable to hope that the bubbling process will ended in
finite step. But we need to be a little careful, as in Einstein case
\cite{Nakajima92}, there could be some overlap if any singular point
lies in a ball centered at other singular point.

So we next consider the ration $\frac{r_i'}{r_i}$.

Case(i): there exists a subsequence(which we continue to index with
$i$) satisfying $r_i'<Cr_i$ for some constant $C$.

Case(ii):
\[\lim_{i\rightarrow\infty}\frac{r_i'}{r_i}=\infty.\]

In Case(i) we proceed as follows: We claim that for $i$ sufficiently
large, the balls $B(x_i,2r_i)$ (from the first subsequence) and
$B(x_i',2r_i')$ (from the second) must disjoint because of choice in
(\ref{AVoRaA1}). To see this, if $B(x_i,2r_i)\cap
B(x_i',2r_i')\neq\emptyset$, then $B(x_i',2r_i')\subset
B(x_i,6r_i')$. Then (\ref{UVGLMD}) and (\ref{AVoRaA1}) implies that
\begin{eqnarray*}
 2600A_1(r_i')^4&=&\textrm{Vol}(B(x_i',r_i'))\\
 &<&\textrm{Vol}(B(x_i',2r_i'))<Vol(B(x_i,6r_i')) \\
  &\leq&2A_1(6(r_i')^4)=2592(r_i)^4
\end{eqnarray*}
which is a contradiction (note the last inequality is true for $i$
sufficiently large since the volume is continuous under
$C^{1,\alpha}$ topology (even with finite singularity), which is
valid only in the Case (i)).

In case (ii), if the balls $B(x_i,2r_i)$ (from the first
subsequence) and $B(x_i',2r_i')$ (from the second) are disjoint for
all $i$ sufficiently large, then it will account at leat
$2\epsilon_0$ of $L^2$ - curvature. Otherwise, we look again at the
scaling so that $r_i'=1: \tilde{g}_i=(r_i')^{-2}g_i$, and basepoint
$x_i'$. Then in this rescaled metric,
\[\textrm{Vol}(B(x_i',1))=2600A_1.\]
As above, we have a limiting smooth multi - fold
$(M_\infty',g_\infty',x_\infty')$, satisfying
\[\textrm{Vol}(B(x_\infty',1))=2600A_1.\]
Since the metric is anti - self - dual with constant curvature, by
the choice of $A_1$, we concluded that
\[\int_{B_{g_\infty'}(x_\infty,2)}|Rm|^2>\epsilon_0^2.\]
There is now a singular point of convergence corresponding  to the
balls $B(x_i,r_i)$ in the first subsequence. But since we are in
Case (ii) with $\lim_{i\rightarrow\infty}\frac{r_i'}{r_i}=\infty$,
in the $g_i'$ metric, these balls must limit a point in $M_\infty'$.
The only possibility is that the original sequence satisfies
\[\int_{B_{g_i}(x_i',2r_i')}|Rm(g_i)|^2>2\epsilon_0^2,\]
for all $i$ sufficiently large.

We repeat the above procedure, considering possible Cases (i) and
(ii) at each step. At the $k\,\textrm{th}$ step, we can always
account for at least $k\,\epsilon_0$ of $L^2-$ curvature. The
process must terminate in finitely many steps from the bound
$\|Rm(g_i)\|_{L^2}<\Lambda_1$. This contradicts (\ref{CUVAC}), which
finishes the proof. \qed

We note that, it may happen that $(M_\infty,g_\infty)$ is a smooth
Riemannian manifold, but the convergence is not in the
$C^{1,\alpha}$ topology. In fact, the curvature concentrate part,
corresponding to some nontrivial $2$ - cycles in $M$,  maybe shrink
off in the limit. The topology `decrease' and the singularity take
away a certain quantity of energy of curvature.

\proposition Let $(M, g_i)$ satisfies hypothesis in theorem
\ref{MTCC}, then
\[\lim_{i\rightarrow \infty}\int_M|Rm(g_i)|\geq \int_{\mathcal{R}\subset M_\infty}|Rm(g_\infty)|^2,\]
with inequality if only if $M_\infty$ is $C^{1,\alpha}$ manifold
diffeomorphic to $M$, and the convergence is in the $C^{1,\alpha}$
Cheeger - Gromov topology.

\proof It is a straightforward consequence of the Bubble analysis.
Since the converge is taking in the $C^{1,\alpha}\cap W^{2,p}$, then
measure $|Rm(g_i)|^2dv_{g_j}$ converge to
\[|Rm(g_\infty)|dv_{g_\infty}+\sum_{x_i\in\mathcal{S}_\infty} a_i\delta_{x_i}\]
in the $L^p$ sense, where $\delta_{x_i}$ is the Dirac measure
supported at $x_1$, and $a_k$ is given by
\[a_k=\sum_{(N_k^*,h_k^*)}\int_{(N_k^*,h_k^*)}|Rm(h_k^*)|^2\]
Here $(N_k^*,h_k^*)$ is the bubble tree associated to the singular
point $x_k$, see clear description of bubble tree in
\cite{Nakajima92}. The equality implies there are no curvature
concentration occurred and thus no singularities in the limit. \qed

\question Does the maximum of sectional curvature (curvature
singularity) always occurred in the singular part $\mathcal
{S}_i(r)$ in theorem \ref{CoLiTh}?

Since there is no  $\epsilon -$ regularity theorem, in appearance,
curvature may blow up even if small energy. If this phenomena does
not occur, combined with the results and methods in
\cite{Anderson89}, one may show that under certain circumstances, in
fact orbifold singularities do not arise in the limit.

For instance,
%
in dimension $4$, one can introduce the bound
\begin{equation}\label{IAESB}
    \inf\{\textrm{area}(\Sigma): [\Sigma]\neq 0 \in H_2(M,\mathbb{Z})/\textrm{torsion}\geq a>0.\}
\end{equation}
Then the deepest bubble $N$ associated to each singularity will be a
non flat ALE space with $\textrm{dim} H_2(N)\neq 0$. Then there
exists $[\Sigma_i]\neq 0$ in $H_2(M_i;\mathbb{Z})$ with area
$\textrm{Area}(\Sigma_i)\rightarrow 0$ as $i\rightarrow \infty$,
i.e. the metric to become degenerate along $[\Sigma_i]$,  which
contradict to (\ref{IAESB}). Consequently, the convergence of
Riemannian manifolds in theorem \ref{MTCC} will be in the
$C^{1,\alpha}$ topology.

\textsc{Yiyan Xu, School of Mathematical Sciences, Peking
University, Beijing, China, 100871}
\hfill\href{mailto:xuyiyan@math.pku.edu.cn}{\textcolor[rgb]{0.00,0.00,1.00}{
xuyiyan@math.pku.edu.cn}}

 \end{document}